\documentclass[12pt]{article}
\usepackage{amsmath, amsthm, amssymb, amsfonts, enumerate,bm}
\usepackage[figuresright]{rotating}
\usepackage{natbib}

\begin{document}
\newcommand{\abs}[1]{\left\vert#1\right\vert}
\newcommand{\set}[1]{\left\{#1\right\}}
\newcommand{\eps}{\varepsilon}
\newcommand{\To}{\rightarrow}
\newcommand{\inv}{^{-1}}
\newcommand{\ihat}{\hat{\imath}}
\newcommand{\var}{\mbox{Var}}
\newcommand{\sd}{\mbox{SD}}
\newcommand{\cov}{\mbox{Cov}}
\newcommand{\f}{\frac}
\newcommand{\fI}[1]{\frac{1}{#1}}
\newcommand{\what}[1]{\widehat{#1}}
\newcommand{\hhat}[1]{\what{\what{#1}}}
\newcommand{\wtilde}[1]{\widetilde{#1}}
\newcommand{\bdot}{\bm{\cdot}}
\newcommand{\Th}{\theta}
\newcommand{\qmq}[1]{\quad\mbox{#1}\quad}
\newcommand{\qm}[1]{\quad\mbox{#1}}
\newcommand{\tr}{\mbox{tr}}
\newcommand{\logit}{\mbox{logit}}
\newcommand{\noi}{\noindent}
\newcommand{\bni}{\bigskip\noindent}
\newcommand{\bul}{$\bullet$ }
\newcommand{\bias}{\mbox{bias}}
\newcommand{\conv}{\mbox{conv}}
\newcommand{\spn}{\mbox{span}}
\newcommand{\colspace}{\mbox{colspace}}
\newcommand{\mX}{\mathcal{X}}
\newcommand{\bbR}{\mathbb{R}}

\newtheorem{theorem}{Theorem}[section]
\newtheorem{corollary}{Corollary}[section]
\newtheorem{conjecture}{Conjecture}[section]
\newtheorem{proposition}{Proposition}[section]
\newtheorem{lemma}{Lemma}[section]
\newtheorem{definition}{Definition}[section]
\newtheorem{example}{Example}[section]
\newtheorem{remark}{Remark}[section]

\title{{\bf\Large A New Characterization of Elfving's Method for High Dimensional Computation}}


\author{Jay Bartroff\\
\small{Department of Mathematics, University of Southern California, Los Angeles, California, USA}\\ \small{\textsf{bartroff@usc.edu}}}   

\date{}
\maketitle

\abstract{We give a new characterization of Elfving's \citeyearpar{Elfving52} method for computing $c$-optimal designs in $k$ dimensions which gives explicit formulae for the $k$ unknown optimal weights and $k$ unknown signs in Elfving's characterization.  This eliminates the need to search over these parameters to compute $c$-optimal designs, and thus reduces the computational burden from solving a family of optimization problems to solving a single optimization problem for the optimal finite support set.  We give two illustrative examples: a high dimensional polynomial regression model and a logistic regression model, the latter showing that the method can be used for locally optimal designs in nonlinear models as well.}

\section{Introduction}

\citet{Elfving52} gave an elegant geometric characterization of $c$-optimal designs for linear regression models by associating the design's support points and weights with a weighted and signed linear combination of those points, for a certain optimal choice of signs \citep[see also][]{Chernoff72,Chernoff99}. This can be used to immediately find $c$-optimal designs in 1 or 2 dimensions, but the state of the art for $k\ge 3$ dimensions is an algorithm of \citet{Lopez-Fidalgo04} that involves searching over the unknown signs and support points.  After introducing our notation and assumptions in Section~\ref{sec:set-up}, in Section~\ref{sec:result} we give a new result on Elfving's method  which gives explicit formulae for the optimal signs and weights for any candidate set of support points, thus reducing the computational cost of the design, greatly when $k$ is large. Torsney and coauthors \citep{Torsney81,Kitsos88,Pukelsheim91} have given a different characterization of the optimal weights, but the novelty of the current method is (i) by giving both the optimal weights and signs, the current method allows fast computation by using Elfving's geometric method along the lines of \citet{Lopez-Fidalgo04} without needing to search over all signs, and (ii) the formulae for the weights and signs also hold for suboptimal support points, which are often used in practice.
In Section~\ref{sec:examples} we discuss two examples: high-dimensional polynomial regression and estimating the turning point of a quadratic logistic regression model, in which we use Ford,  Torsney, and Wu's \citeyearpar{Ford92} results for locally $c$-optimal designs in generalized linear models.

In addition to the references already mentioned, there has been much recent research on Elfving's problem, of which \citet{Studden05} gives a review. \citet{Dette93} showed that the popular $D$-optimal designs can be written as a sequence of $c$-optimal design problems and also considered a Bayesian approach. Later, \citet{Dette97} proposed ``standardized'' optimality criteria to optimize estimation of a set of parameter values, in various senses, which involve $c$-optimal designs as their constituents. Elfving's method has been generalized by Dette and coauthors in several other directions as well in \citet{Dette93b}, \citet{Dette93d, Dette93c, Dette96}, \citet{Dette08}, \citet{Dette09}, and \citet{Dette10}.

\citet{Fan03} characterized a certain generalized inverse of a candidate $c$-optimal design that can be used to assess optimality. In a recent application,  \citet{Bartroff10b,Bartroff10d} have proposed using $c$-optimal designs as ``base designs'' in dynamic programming-based dose-finding clinical trial designs which are computed repeatedly in sequential Monte Carlo simulation, putting a premium on efficient computation.

\subsection{Notation and assumptions}\label{sec:set-up}
We assume that the scalar-valued response variable $Y$ is normally distributed with mean~$\theta'x$ and variance $\sigma^2$, where prime denotes transpose, $\theta\in\mathbb{R}^k$ is the unknown parameter, and $x$ is the design point taking values in the compact design space $\mathcal{X}\subseteq \mathbb{R}^k$. Our main assumption is that any collection of linearly independent~(LI) vectors $x_1,\ldots,x_\ell\in\mathcal{X}$, $\ell<k$, can be completed to a LI set $x_1,\ldots,x_\ell,x_{\ell+1},\ldots,x_k\in\mathcal{X}$ of size $k$; this is not a strong assumption and any $\mathcal{X}$  not satisfying this can be reduced to a smaller dimensional case, in some sense. We define the information matrix of a given measure $\xi$ on $\mathcal{X}$ as
\begin{equation*}
M(\xi)=\int_\mX  xx'\xi(dx).
\end{equation*} For a given nonzero vector $c\in\mathbb{R}^k$, the $c$-optimal design  problem is to find a design measure $\xi$ achieving
\begin{equation}\label{uncon-problem}
\min_{\xi\in\Xi} \Psi(\xi),\qmq{where} \Psi(\xi)=c'M(\xi)^{-1}c.
\end{equation} 
Here  $\Xi$ is the class of designs for which $c'\theta$ is \emph{estimable}, i.e., 
\begin{equation*}
\Xi=\left\{\xi: \int_{\mathcal{X}}\xi(dx)=1\qmq{and} c=M(\xi)d \qm{for some $d\in\mathbb{R}^k$}\right\}.
\end{equation*}
The class $\Xi$ may admit singular designs in which case $M(\xi)^{-1}$ in \eqref{uncon-problem} means a generalized inverse \citep[see][p.~469]{Seber03} and the quadratic form $c' M(\xi)^{-1}c$ is invariant under which generalized inverse is chosen. Also define
\begin{equation*}
\Xi(x_1,\ldots,x_\ell)=\{\xi\in\Xi: \mbox{supp}(\xi)\subseteq \{x_1,\ldots,x_\ell\}\}.
\end{equation*}

\citet{Elfving52} gave the following elegant solution to the problem (\ref{uncon-problem}). Let $\mathcal{X}^-$ denote the reflection of $\mathcal{X}$ through the origin and $\mathcal{E}=\mbox{conv}(\mathcal{X}\cup \mathcal{X}^-)$, the convex hull of $\mathcal{X}\cup \mathcal{X}^-$.   Points in $\mathcal{E}$, the so-called Elfving set, can be written in the form $\sum_{i=1}^{n} p_i\eps_ix_i$, where $p_i\ge 0$, $\sum_{i=1}^n p_i=1$, $\eps_i\in\{-1,+1\}$, and $x_i\in \mathcal{X}$. Let $z=\sum_{i=1}^{n} p_i\eps_ix_i\in \mathcal{E}$ be the point on the ray 
\begin{equation}\label{Rc}
R(c)=\{\gamma c\in\mathbb{R}^k: \gamma\in\mathbb{R}\}
\end{equation} that is furthest from the origin. Then the design $\{(x_i,p_i)\}_{i=1}^n$, which puts weight $p_i$ on point $x_i$, $i=1,\ldots,n$, is $c$-optimal. Moreover, $n$ can be taken to be $k$ by Carath\'{e}odory's theorem, and
\begin{equation}\label{ElfPsic/z}
\Psi(\{(x_i,p_i)\})=||c||^2/||z||^2.\end{equation} See \citet{Chernoff72,Chernoff99}.

\section{A new characterization of Elfving's method and an improved algorithm}\label{sec:result}
For LI vectors $x_1,\ldots,x_\ell\in\mathcal{X}$, $\ell\le k$, such that $c\in\spn(x_1,\ldots,x_\ell)$, define the design $\xi^*(x_1,\ldots,x_\ell)$ and the point $z^*(x_1,\ldots,x_\ell)\in \mathbb{R}^k$ as follows:
\begin{align}
\xi^*(x_1,\ldots,x_\ell)&=\{(x_i,p_i^*)\}_{i=1}^\ell\label{xi*}\\
z^*(x_1,\ldots,x_\ell)&=\sum_{i=1}^\ell \eps_i^*p_i^*x_i,\label{z*}
\end{align} where
\begin{align}
\eps_i^*&=\eps_i^*(x_1,\ldots,x_\ell)=\mbox{sign}\left(x_i'M^{-1}c\right),\quad i=1,\ldots,\ell,\label{epsi}\\
M&=M\left(\{(x_i,1/\ell)\}_{i=1}^\ell \right)\;\mbox{and $M^{-1}$ any generalized inverse,}\label{M}\\ 
p_i^*&=p_i^*(x_1,\ldots,x_\ell)=\abs{x_i'M^{-1}c}\left/\sum_{j=1}^\ell \abs{x_j'M^{-1}c}\right.,\quad i=1,\ldots,\ell.\label{pi}
\end{align}
Our main result, Theorem~\ref{thm:phi}, is that solving the single optimization problem of maximizing the function
\begin{equation*}
\phi(x_1,\ldots,x_k)=||z^*(x_1,\ldots,x_k)||^2
\end{equation*}
over $\mathcal{S}=\{(x_1,\ldots,x_k)\in\mathcal{X}^k: x_1,\ldots,x_k\;\mbox{are LI}\}$ gives  a $c$-optimal design. The proof of Theorem~\ref{thm:phi} is given in the Appendix, as well as proofs of Proposition~\ref{prop:xi*}, which establishes the properties of $\xi^*(x_1,\ldots,x_\ell)$ and $z^*(x_1,\ldots,x_\ell)$, and Lemma~\ref{lem:span=C} relating $\mbox{span}(x_1,\ldots,x_\ell)$ to the column space of $M(\{(x_i,q_i)\}_{i=1}^\ell)$.

\begin{theorem}\label{thm:phi}
Let $(x_1^*,\ldots,x_k^*)=\arg\max_\mathcal{S} \phi(x_1,\ldots,x_k)$.  Then the design $$\xi^*(x_1^*,\ldots,x_k^*)=\{(x_i^*,p_i^*)\}_{i=1}^k$$ is $c$-optimal, where $p_i^*=p_i^*(x_1^*,\ldots,x_k^*)$ are given by \eqref{pi}.
\end{theorem}

\bigskip

\noindent\textbf{Remarks.} 
\begin{enumerate}
\item There is nothing special about the weights $1/\ell$ in \eqref{M}. In fact, any strictly positive constants $a_1,\ldots,a_\ell$ (not necessarily summing to $1$) can be used if we replace $M$ in \eqref{epsi}-\eqref{M} by $\sum_{i=1}^\ell a_ix_ix_i'$ and $p_i^*$ in \eqref{pi} by
\begin{equation}
p_i^*=a_i\abs{x_i'M^{-1}c}\left/\sum_{j=1}^\ell a_j\abs{x_j'M^{-1}c}\right..
\end{equation} 
\item In practice, one is often constrained to apply a suboptimal design with support on some given set.  Thus it is convenient to be able to compute the optimal weights for the given support set. Claim~\ref{Psimin} of Proposition~\ref{prop:xi*}, given in the Appendix, shows that \eqref{pi} can be used to find these.
\item For many models of interest, the vectors $x\in\mathcal{X}\subseteq\mathbb{R}^k$ are actually functions $x=x(u)$ of a scalar $u\in\mathcal{U}\subseteq\mathbb{R}$, for some closed interval $\mathcal{U}$.  In this case the optimization problem of size $k^2$ (i.e., $k$ vectors in $\mathbb{R}^k$) in Theorem~\ref{thm:phi} reduces to an optimization problem of size~$k$.
\end{enumerate}

\section{Illustrative examples}\label{sec:examples}
\subsection{Polynomial regression}
The method of Theorem~\ref{thm:phi} was used to compute the $c$-optimal designs in the polynomial model $Y\sim\mathcal{N}(\sum_{i=1}^{k}\theta_iu^{i-1},\sigma^2)$ for $k=6,\ldots,10$, $u\in[-1,1]$, and for $c$ equal to the $k$ standard basis vectors $e_1=(1,0,\ldots,0)',\ldots,e_k=(0,\ldots,0,1)'$.  The corresponding designs for $k=3,4$, and $5$ appear in Lopez-F\'idalgo and Rodriguez-Diaz\footnote{We have verified the values in \citet[Table~5]{Lopez-Fidalgo04} using our method, except for the value of $\Psi(\xi^*)$ in the case $c=e_2$, $k=3$ ($d=2$ in their notation), which we calculate to be 1 rather than the value~4 given there.} \citeyearpar[Table~5]{Lopez-Fidalgo04}, and the $k=1,2$ cases are well known. Here, in Table~\ref{table:poly}, and in the next section, we abuse the above notation slightly by parametrizing designs $\xi^*$ by $u\in[-1,1]$ rather than $x=x(u)=(1,u,u^2,\ldots,u^{k-1})'\in\bbR^k$. A few interesting properties of these designs are (i) all are symmetric about $u=0$, (ii) for all $k$, the $e_1$-optimal design is $\xi^*=\{(0,1)\}$ with $\Psi(\xi^*)=1$, and (iii) for all $k$, the $e_j$-optimal design has positive weight at $u=0$ if and only if $j$ is odd. These and other properties were noted by \citet{Studden68} who solved this problem explicitly using Tchebycheff systems.

The algorithm was implemented in \textsf{Matlab} using its \textsf{fminsearch} optimizer. On a 2.6~GHz laptop computer, the design in Table~\ref{table:poly} that took the longest to compute took roughly 3 seconds. We have computed these designs up to dimension $k=20$ using the same method with the longest taking roughly 10 seconds.

\subsection{Quadratic logistic regression}\label{sec:logistic}

We consider $c$-optimal designs for finding the turning point of a quadratic logistic regression model. Suppose $Y$ takes the values 0 or 1 according to 
\begin{equation}\label{logit}
P_\theta(Y=1|u)=1/(1+e^{-(\theta_1+\theta_2 u+\theta_3 u^2)}),\end{equation} where $\theta=(\theta_1,\theta_2,\theta_3)'$ and the scalar $u$ takes values in some interval, taken here to be $[-1,1]$. For example, this quadratic logistic regression is popular for modeling the response to treatments for diseases other than cancer, for which efficacy is assumed to have a unimodal relationship to dose~$u$ \citep[e.g.,][p.~686]{Thall04}. In settings such as this it is assumed that $\theta_3<0$ and it is desired to estimate the dose $\eta=-\theta_2/(2\theta_3)$ maximizing probability of response $Y=1$. Given an estimate $\what{\theta}$ of $\theta$ with $\what{\theta}_3\ne 0$, we consider this as a locally (i.e., depending on $\theta$) $c$-optimal design problem by first writing
\begin{equation*}
\eta\approx -\what{\theta}_2/2\what{\theta}_3-(\theta_2-\what{\theta}_2)/(2\what{\theta}_3)+(\theta_3-\what{\theta}_3)\what{\theta}_2/(2\what{\theta}_3^2),
\end{equation*} and hence taking $c$ to be $(0,-1/(2\what{\theta}_3),\what{\theta}_2/(2\what{\theta}_3^2))'$, or equivalently, 
\begin{equation}\label{cforeta}
c=(0,-\what{\theta}_3,\what{\theta}_2)'.\end{equation} \citet{Chaloner89} considered a Bayesian approach to this problem in the linear model.

Although the method described above in Section~\ref{sec:result} is for the linear model described in Section~\ref{sec:set-up}, it can be applied to find locally $c$-optimal designs in the nonlinear model~\eqref{logit}, and other generalized linear models, without further extension by applying Ford et al.'s \citeyearpar{Ford92} results which show that, in this case, the locally $c$-optimal design can be found by applying Elfving's method to a transformed version of the problem as if it were a linear model, described next. For a more general description of their method, which covers other models and optimality criteria, we refer readers to \citet{Ford92}.

In the above notation, $x=x(u)=(1,u,u^2)'$ and $\mX=\{x(u): u\in[-1,1]\}$. Given $\what{\theta}$, let  $B=B(\what{\theta})=||\what{\theta}||U$ where $U=U(\what{\theta})$ is a $3\times 3$ orthonormal matrix with third row equal to $\what{\theta}'/||\what{\theta}||$. Then the locally $c$-optimal design can be found by finding the $(Bc)$-optimal design in the transformed design space
\begin{equation}\label{G}
\mathcal{G}=\left\{w(z_3) z: z=Bx, x\in\mX\right\}\qmq{where} w(\zeta)=e^{\zeta/2}/(1+e^\zeta),\end{equation} and then transforming back to $\mX$ via $B^{-1}$. We give an example of this transformation and calculation.

\begin{figure}[p]
\begin{center}
\scalebox{.7}{\includegraphics{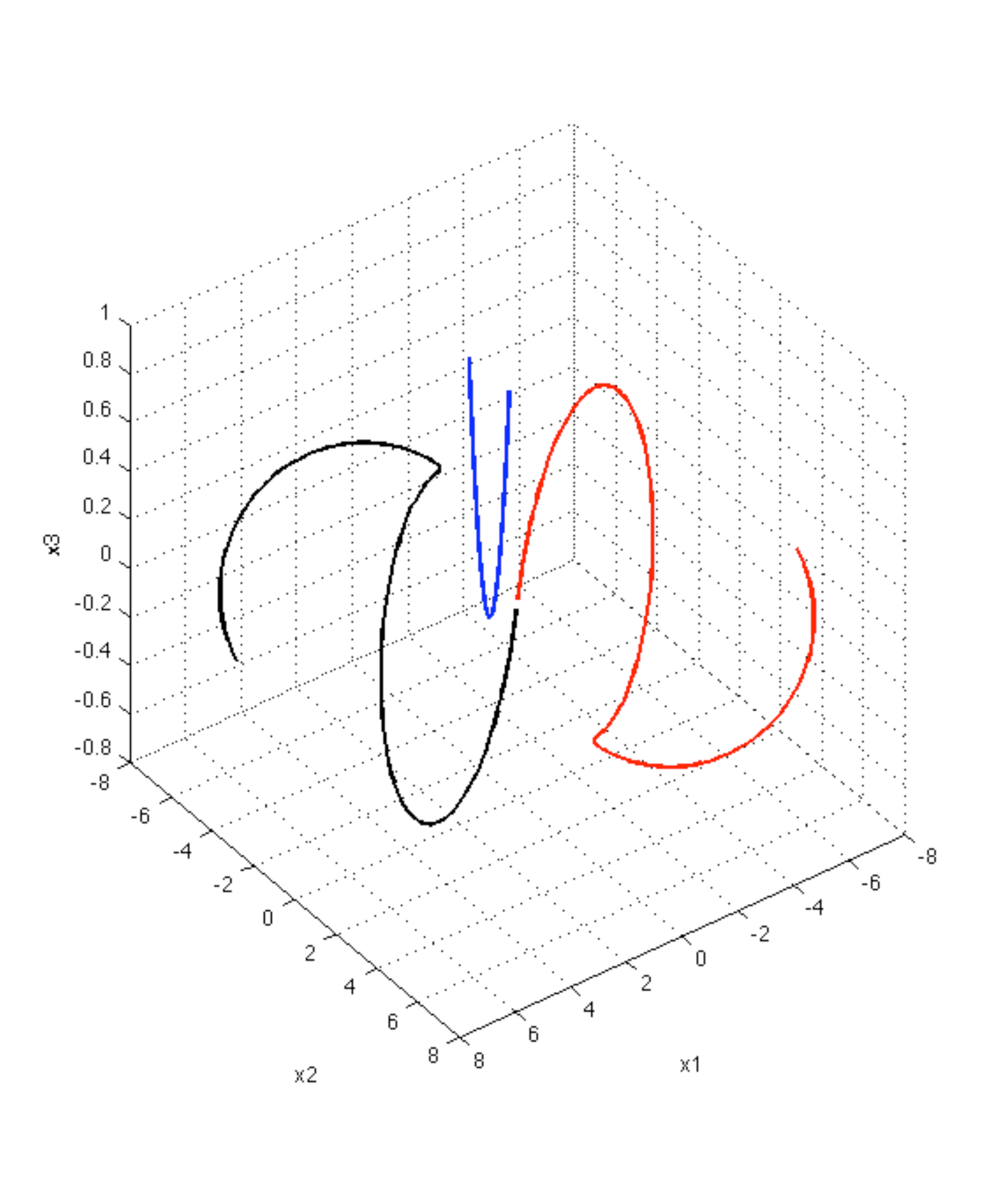}}
\caption{For the quadratic logistic regression example in Section~\ref{sec:logistic}, $\mX$ (blue), $\mathcal{G}$ (red) given by \eqref{G} with $B$ equal to \eqref{B}, and $\mathcal{G}^-$ (black).}
\label{fig:X&Xtild}
\end{center}
\end{figure}

\begin{figure}[p]
\begin{center}
\scalebox{.9}{\includegraphics{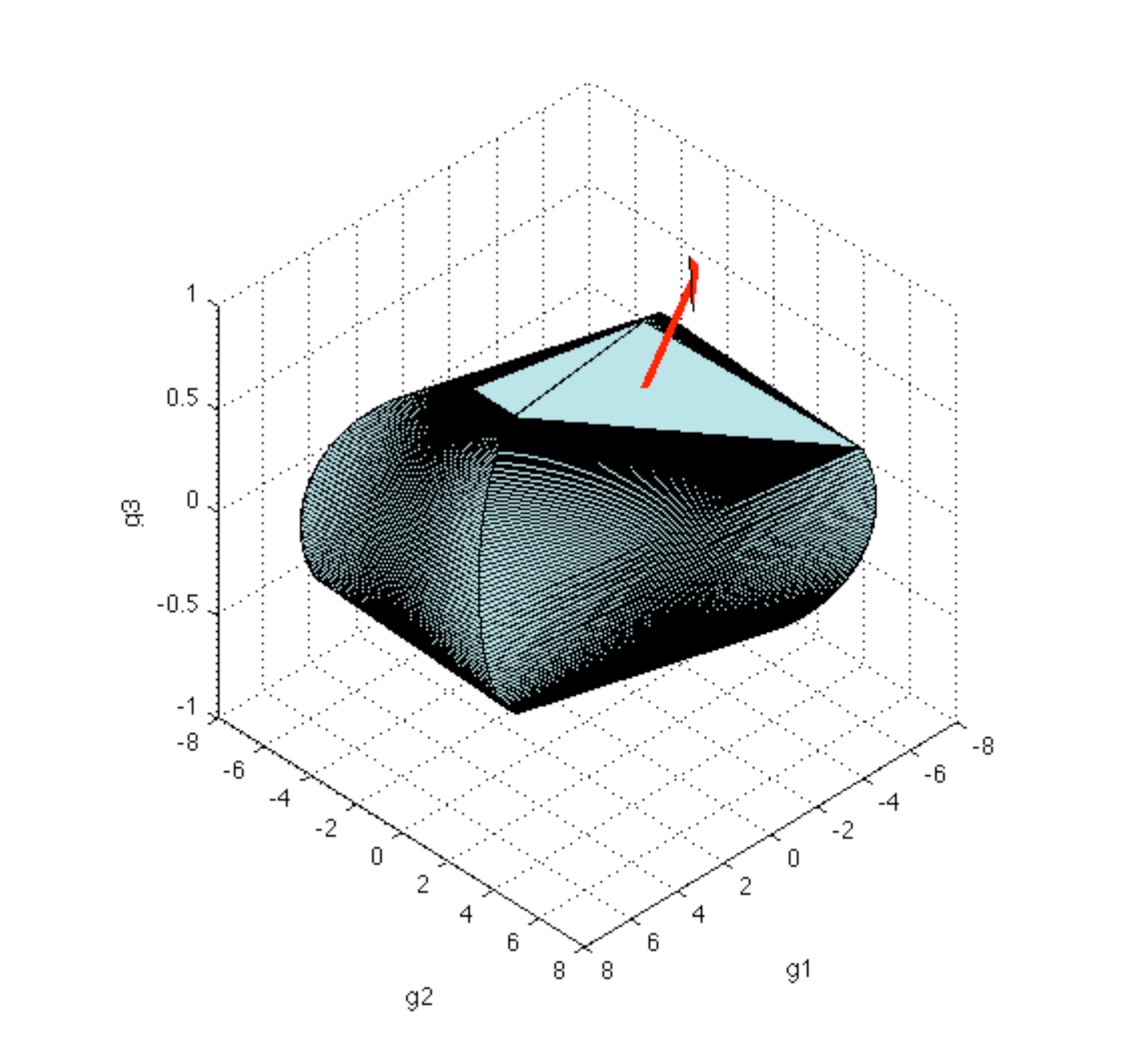}}
\caption{$\mbox{conv}(\mathcal{G}\cup \mathcal{G}^-)$ for the quadratic logistic regression example in Section~\ref{sec:logistic}, with $\mathcal{G}$ given by \eqref{G} with $B$ equal to \eqref{B} and with red arrow $Bc=(-2.403,2.398,1.197)'$.}
\label{fig:convw/c}
\end{center}
\end{figure}

Letting $\what{\theta}=(2,-6,-9)'$, a $B$ satisfying the above assumptions is
\begin{equation}\label{B}
B=\left[
\begin{array}{ccc}
 10.816 &1.110&1.664   \\
  0&9.153&-6.102\\
  2&-6&-9 
\end{array}
\right],
\end{equation} to 3 decimal places. Figure~\ref{fig:X&Xtild} shows $\mX$, $\mathcal{G}$ for this $B$, and the reflection $\mathcal{G}^-$ of $\mathcal{G}$ through the origin, and Figure~\ref{fig:convw/c} shows $\mbox{conv}(\mathcal{G}\cup \mathcal{G}^-)$. For the sake of illustration we choose $c=(-.195,.1,-.243)'$, although this does not have first component equal to zero as in \eqref{cforeta}, and this gives $Bc=(-2.403,2.398,1.197)'$, the red arrow in Figure~\ref{fig:convw/c}. The method is then used to quickly compute the $(Bc)$-optimal design on $\mathcal{G}$ and transforming back to $\mX$ gives, parametrizing by $u$, $$\xi^*=\{(-1,.135), (.181,.194), (.452,.671)\}.$$ From Figure~\ref{fig:convw/c} it is clear that the $(Bc)$-optimal design on $\mathcal{G}$ will be supported on the point $\overline{g}$ at the end of the red curve in Figure~\ref{fig:X&Xtild} with positive second coordinate, plus two interior points on the curve.  Noting that $\overline{g}$ corresponds to the $u=-1$ end of the curve $\mX$ under the transformation $\mX\To\mathcal{G}$, it is not surprising then that $\xi^*$ is supported on $u=-1$ plus two interior points in $(-1,1)$.

\section{Discussion}

We have given a new characterization of Elfving's method which affords more efficient computation, particularly in high dimensions where visualization of the Elfving set $\mathcal{E}$ is difficult.  Dette's \citeyearpar{Dette93} extension of Elving's method to the popular $D$-optimality criterion holds the possibility extending the current computational method to find $D$-optimal designs. The approach of \citet{Ford92} to generalized linear models applied in Section~\ref{sec:logistic} would apply there as well since the $D$-optimality criterion is also $B$-invariant. 

Another tantalizing area of extension is constrained $c$-optimal designs \citep{Cook94,Cook95}. With linear constraints on the design measure, as in \citet[][Section~1.3]{Cook95}, it can be shown that the constrained problem has similar ingredients to the standard $c$-optimal design problem. For example, the subset of points in $\mathcal{E}$ which satisfy the constraint is convex and symmetric about the origin, so one may suspect that, like the standard problem, the point where the ray $R(c)$ pierces this subset is a natural candidate for the constrained $c$-optimal design. However, difficulties remain in this approach, such as showing that an Elfving-like characterization like \eqref{ElfPsic/z} holds for the constrained problem. \citet{Stigler71} has pointed out similar difficulties for the constrained $D$-optimal problem.

\section*{Acknowledgements} The author thanks two reviewers for their helpful comments.  This work was partially supported by grants from the National Science Foundation (Number DMS-0907241) and the National Security Agency (Number H98230-11-1-0162).

\section*{Appendix: Properties of the design \eqref{xi*} and proof of Theorem~\ref{thm:phi}}

Before giving the proof of Theorem~\ref{thm:phi}, we state and prove Proposition~\ref{prop:xi*} which gives the properties of the design $\xi^*$, given by \eqref{xi*}, and the ``Elfving point'' $z^*$, given by \eqref{z*}, needed to prove Theorem~\ref{thm:phi}.  We state this as a separate proposition because it may be of interest on its own. In particular, Claim~\ref{Psimin} of Proposition~\ref{prop:xi*} shows that the weights \eqref{pi} of the design $\xi^*(x_1,\ldots,x_\ell)$ are optimal for the given support set $x_1,\ldots,x_\ell$, showing that the formulae \eqref{M}-\eqref{pi} can be used to compute sub-optimal designs with prescribed support.

\begin{proposition}\label{prop:xi*}
Let $x_1,\ldots,x_\ell\in\mathcal{X}$, $\ell\le k$, be LI such that $c\in\mbox{span}(x_1,\ldots,x_\ell)$. Then the following hold.
\begin{enumerate}
\item\label{invariant} $\xi^*(x_1,\ldots,x_\ell)$ and $z^*(x_1,\ldots,x_\ell)$ are invariant with respect to the choice of generalized inverse in \eqref{M}-\eqref{pi}
\item $z^*(x_1,\ldots,x_\ell)\in R(c)$\label{zinR}
\item $\xi^*(x_1,\ldots,x_\ell)\in \Xi$\label{xiinXi}
\item$\Psi(\xi^*(x_1,\ldots,x_\ell))=||c||^2/||z^*(x_1,\ldots,x_\ell) ||^2$\label{Psi=c/z^2}
\item $\xi^*(x_1,\ldots,x_\ell)=\arg\min_{\xi\in\Xi(x_1,\ldots,x_\ell)} \Psi(\xi)$ \label{Psimin}
\end{enumerate}
\end{proposition}

\noindent\textbf{Proof.}  For any design $\xi=\{(x_i,p_i)\}_{i=1}^\ell\in\Xi$, we have that $c=M(\xi)d$ for some $d$, hence
\begin{multline}\label{c=Md}
c=M(\xi)d=M(\xi)M(\xi)^{-1}M(\xi)d=M(\xi)M(\xi)^{-1}c=\left(\sum_{i=1}^\ell p_ix_ix_i'\right)M(\xi)^{-1}c\\
=\sum_{i=1}^\ell p_i(x_i'M(\xi)^{-1}c)x_i.
\end{multline}
Letting $M$ be as in \eqref{M}, we have that $c\in \mbox{span}(x_1,\ldots,x_\ell)\subseteq\mbox{colspace}(M)$ by Lemma~\ref{lem:span=C}, which is further below in this appendix. Let $\gamma=\ell/\sum_{j=1}^\ell\abs{x_j'M^{-1}c}$ and denote $z^*(x_1,\ldots,x_\ell)$ simply by $z^*$. Applying \eqref{c=Md} with $\xi=\{(x_i,1/\ell)\}_{i=1}^\ell$ gives 
\begin{equation}\label{c=M*d}
c=\ell^{-1}\sum_{i=1}^\ell(x_i'M^{-1}c)x_i=\gamma^{-1}\sum_{i=1}^\ell \eps_i^*p_i^*x_i=\gamma^{-1}z^*.
\end{equation}
This gives Claim~\ref{zinR}, and also Claim~\ref{xiinXi}, as follows.  By relabeling if necessary, let $x_1,\ldots,x_m$, $m\le \ell$, be the $x_1,\ldots,x_\ell$ corresponding to nonzero $p_i^*$'s.  Then \eqref{c=M*d} shows that $c\in\spn(x_1,\ldots,x_m)\subseteq\colspace(M(\xi))$ by Lemma~\ref{lem:span=C}. 

Denote $\xi^*(x_i,\ldots,x_\ell)=\{(x_i,p_i^*)\}_{i=1}^m$ simply by $\xi^*$. For Claim~\ref{Psi=c/z^2}, apply \eqref{c=Md} with $\xi=\xi^*$ to get
$$c=\sum_{i=1}^mp_i^*(x_i'M(\xi^*)^{-1}c)x_i,$$
and subtracting this from the second-to-last expression in \eqref{c=M*d} (with $m$ in place of $\ell$) gives
\begin{equation*}
\sum_{i=1}^m\left(\eps_i^*/\gamma-x_i'M(\xi^*)^{-1}c\right)p_i^*x_i=0.
\end{equation*} Since the $x_i$ are LI and all the $p_i^*$ in the sum are positive, this gives
\begin{equation}
x_i'M(\xi^*)^{-1}c=\eps_i^*/\gamma\quad\mbox{for all $i=1,\ldots,m$.}
\end{equation} Then
\begin{multline*}
\Psi(\xi^*)=c'M(\xi^*)^{-1}c=\left(\gamma^{-1}\sum_{i=1}^m \eps_i^* p_i^*x_i\right)'M(\xi^*)^{-1}c= \gamma^{-1}\sum_{i=1}^m \eps_i^*p_i^*(x_i'M(\xi^*)^{-1}c) \\= \gamma^{-1}\sum_{i=1}^m \eps_i^*p_i^*(\eps_i^*/\gamma)
=\gamma^{-2}=||c||^2/||z^*||^2,
\end{multline*} this last by \eqref{c=M*d}, proving Claim~\ref{Psi=c/z^2}.

To show Claim~\ref{Psimin}, subtract \eqref{c=Md} from the second-to-last expression in \eqref{c=M*d} to get
\begin{equation*}
\sum_{i=1}^\ell\left[\eps_i^*p_i^*/\gamma-p_i(x_i'M(\xi)^{-1}c)\right]x_i=0
\end{equation*} which, by linear independence, gives
\begin{equation}\label{qipi}
p_i(x_i'M(\xi)^{-1}c)=\eps_i^*p_i^*/\gamma\quad\mbox{for all $i=1,\ldots,\ell$.}
\end{equation}
 Then, using \eqref{c=Md} and \eqref{qipi}, 
\begin{multline}
\Psi(\xi)=c'M(\xi)^{-1}c=\left[\sum_{i=1}^\ell p_i(x_i'M(\xi)^{-1}c)x_i\right]'M(\xi)^{-1}c=\sum_{i=1}^\ell p_i(x_i'M(\xi)^{-1}c)^2\\
=\sum_{i:\; p_i>0} p_i(x_i'M(\xi)^{-1}c)^2 =\gamma^{-2}\sum_{i:\; p_i>0} (p_i^*)^2/p_i.\label{Psi=p/q}
\end{multline} Then the method of Lagrange multipliers easily shows that $p_i=p_i^*$ for all $i$ minimizes \eqref{Psi=p/q} subject to $\sum_i p_i=1$, proving Claim~\ref{Psimin}. 

Finally, to prove Claim~\ref{invariant}, suppose a different generalized inverse is used for $M^{-1}$, resulting in $\wtilde{\xi}, \wtilde{\eps}_i, \wtilde{p}_i, \wtilde{\gamma}$, $\wtilde{z}$ instead of $\xi^*,\eps_i^*,p_i^*,\gamma,z^*$. By the same argument leading to \eqref{qipi} we would have 
\begin{equation}\label{qipitilde}
\wtilde{\eps}_i\wtilde{p}_i/\wtilde{\gamma}=\eps_i^*p_i^*/\gamma\quad\mbox{for all $i=1,\ldots,\ell$.}
\end{equation} Also, using \eqref{c=M*d} and Claims~\ref{Psi=c/z^2} and \ref{Psimin},
\begin{equation*}
\gamma=\frac{||z^*||}{||c||}=\left[\Psi(\xi^*)\right]^{-1/2}= \left[\Psi(\xi)\right]^{-1/2}=\frac{||\wtilde{z}||}{||c||}=\wtilde{\gamma},
\end{equation*} hence \eqref{qipitilde} shows that $\eps_i^*p_i^*=\wtilde{\eps}_i\wtilde{p}_i$ for all $i$, and it follows that $\xi^*=\wtilde{\xi}$ and $z^*=\wtilde{z}$.
\qed

\bigskip

\noindent\textbf{Proof of Theorem~\ref{thm:phi}.}  Denote $\xi^*(x_1^*,\ldots,x_k^*)$ simply by $\xi^*$. 
Let $\xi=\{(x_i,p_i)\}_{i=1}^\ell$, $p_i>0$ and $\ell\le k$, be the design given by Elfving's method.
If $x_1,\ldots,x_\ell$ are LI then they can be completed to a LI set $x_1,\ldots,x_k\in\mathcal{X}$ by our assumption in the first paragraph of Section~\ref{sec:set-up}. Then $\xi\in\Xi(x_1,\ldots,x_k)$ and so $\Psi(\xi^*)\le \Psi(\xi)$ by Claim~\ref{Psimin} of Proposition~\ref{prop:xi*}. 

Otherwise, assume that  $x_1,\ldots,x_\ell$ are linearly dependent (LD). In this case we will construct a design $\wtilde{\xi}$ whose support is a strict subset of $x_1,\ldots,x_\ell$ and $\Psi(\wtilde{\xi})\le\Psi(\xi)$. This argument can be repeated until only a LI support set remains, and the previous argument can be applied. By \eqref{ElfPsic/z}, $\Psi(\xi)=||c||^2/||\sum_{i=1}^\ell \eps_ip_ix_i||^2$ for some choice of signs $\eps_i$, and let $z=\sum_{i=1}^\ell \eps_ip_ix_i$. Since $x_1,\ldots,x_\ell$ are LD, $\eps_1x_1,\ldots,\eps_\ell x_\ell$ are also LD  so let $\alpha_1,\ldots,\alpha_\ell$ be constants, not all $0$, such that $\sum_{i=1}^\ell\alpha_i \eps_ix_i=0$. Without loss of generality assume that $\sum_{i=1}^\ell\alpha_i\ge 0$, since otherwise each $\alpha_i$ could be replaced by $-\alpha_i$. Now set $x_0=0$, $\eps_0=1$, $p_0=0$, and $\alpha_0=-\sum_{i=1}^\ell\alpha_i\le 0$. With this, we have that $\alpha_0,\ldots,\alpha_\ell$ are not all $0$ and satisfy
\begin{equation}
\sum_{i=0}^\ell\alpha_i\eps_ix_i= 0\qmq{and} \sum_{i=0}^\ell\alpha_i=0.
\end{equation} (In other words, $\eps_0x_0,\ldots,\eps_\ell x_\ell$ are affinely dependent.) Now set $$\delta=\min\{p_i/\alpha_i: i=0,\ldots,\ell,\;\alpha_i>0\}$$ and $q_i=p_i-\delta\alpha_i$, $i=0,\ldots,\ell$, and note that $q_i\ge 0$ for all $i$ and $$\sum_{i=0}^\ell q_i=\sum_{i=0}^\ell p_i-\delta\sum_{i=0}^\ell \alpha_i=1.$$ Note also that $q_i=0$ for some $0<i\le\ell$ by definition of $\delta$ and by virtue of the fact that $\alpha_0\le 0$. By relabeling, assume that $q_\ell=0$. We have
\begin{equation}\label{z=sumq}
z=z-0=\sum_{i=0}^\ell p_i\eps_i x_i- \delta\sum_{i=0}^\ell \alpha_i\eps_ix_i=\sum_{i=0}^\ell q_i\eps_ix_i =\sum_{i=0}^{\ell-1} q_i\eps_ix_i.
\end{equation} It cannot be that $q_0=1$ since this would imply that $z=0$ by \eqref{z=sumq} and it would follow that $c=0$, contradicting our assumption, hence $q_0<1$. Let $\mX_0$ be any compact set containing $\mX\cup\{0\}$, $$\what{\xi}=\{(x_i,q_i)\}_{i=0}^{\ell-1},\qmq{and}\wtilde{\xi}=\{(x_i,q_i/(1-q_0))\}_{i=1}^{\ell-1}.$$ The former is a design on $\mX_0$ by virtue of its inclusion of $x_0=0$ in its support. Then \eqref{z=sumq} implies that $\Psi(\what{\xi})\le\Psi(\xi)$. Also,
\begin{equation*}
M(\wtilde{\xi})=(1-q_0)^{-1}\sum_{i=1}^{\ell-1}q_ix_ix_i'= (1-q_0)^{-1}\sum_{i=0}^{\ell-1}q_ix_ix_i'=(1-q_0)^{-1}M(\what{\xi})
\end{equation*} so that 
\begin{equation*}
\Psi(\wtilde{\xi})=c'M(\wtilde{\xi})^{-1}c=(1-q_0)\Psi(\what{\xi})\le \Psi(\what{\xi})\le \Psi(\xi),
\end{equation*} thus completing the desired reduction to the smaller support set $x_1,\ldots, x_{\ell-1}$.\qed

\bigskip

\begin{lemma}\label{lem:span=C}
Let $x_1,\ldots,x_\ell\in\mathcal{X}$. Then
 \begin{equation}
\mbox{colspace}\left(M(\{(x_i,q_i)\}_{i=1}^\ell) \right) \subseteq\mbox{span}(x_1,\ldots,x_\ell)\;\mbox{for any $q_i\ge 0$, $\sum_{i=1}^\ell q_i=1$.}
\end{equation}
 Conversely, 
 \begin{equation}
\mbox{span}(x_1,\ldots,x_\ell) \subseteq \mbox{colspace}\left(M(\{(x_i,q_i)\}_{i=1}^\ell) \right)\;\mbox{for any $q_i> 0$, $\sum_{i=1}^\ell q_i=1$.}
\end{equation}
\end{lemma}

\noindent\textbf{Proof.} Let $M= M(\{(x_i,q_i)\}_{i=1}^\ell)$. For any $a\in\mathbb{R}^\ell$, $$Ma=\left(\sum_{i=1}^\ell q_ix_ix_i'\right)a=\sum_{i=1}^\ell \left(q_ix_i' a\right)x_i\in \mbox{span}(x_1,\ldots,x_\ell).$$ Conversely, given constants $a_1,\ldots,a_\ell$ and strictly positive $q_1,\ldots,q_\ell$, we will show that $\sum_{i=1}^\ell a_ix_i \in \mbox{colspace}(M)$. Let $b=(a_1/\sqrt{q_1},\ldots,a_\ell/\sqrt{q_\ell})'\in\mathbb{R}^\ell$. Let $X$ be the $\ell\times k$ matrix with rows $\sqrt{q_1}x_1',\ldots, \sqrt{q_\ell}x_{\ell}'$ and note that $X'X=\sum_{i=1}^\ell q_ix_ix_i'= M$. Let $v$ be the projection of $b$ onto $\mbox{colspace}(X)$ and let $d$ be any vector such that $Xd=v$. Then $(b-v)\perp\mbox{colspace}(X)$, or $X'b=X'v$, so we have
$$\sum_{i=1}^\ell a_ix_i=X'b=X'v=X'Xd=M d\in \mbox{colspace}(M).$$\qed


\def\cprime{$'$}

\begin{table}[p]
\caption{$c$-optimal designs $\xi^*$ for polynomial regression model $Y\sim\mathcal{N}(\sum_{i=1}^{k}\theta_iu^{i-1},\sigma^2)$ and $c$ equal to standard basis vectors $e_1,\ldots,e_k\in\bbR^k$ for $k=6,\ldots,10$.  Parametrizing $\xi^*$ by $u\in[-,1,1]$, column 3 is $\xi^*(0)$ and the pairs $(u,p)$ in columns 4--8 are the positive $u$ values such that $\xi^*(u)=\xi^*(-u)=p$. The support points and weights are given to 3 decimal places.}
\begin{tabular}{c|c|c|lllll|r}\hline\hline
$k$&$j:\; c=e_j$&$\xi^*(0)$&\multicolumn{5}{c|}{$(u,p)$ such that $\xi^*(u)=\xi^*(-u)=p$, $u>0$}&$\Psi(\xi^*)$\\\hline \hline
&1&1&&&&&&1\\
&2&0&(.309, .419)&(.809, .061)&(1, 1/50)&&&25\\
&3&3/8&(.707, 1/4)&(1, 1/16)&&&&64\\
6&4&0&(.309, .265)&(.809, .175)&(1, .06)&&&400\\
&5&1/4&(.707, 1/4)&(1, 1/8)&&&&64\\
&6&0&(.309, 1/5)&(.809, 1/5)&(1, 1/10)&&&256\\
\hline
&1&1&&&&&&1\\
&2&0&(.309, .419)&(.809, .061)&(1, 1/50)&&&25\\
&3&.352&(.5, 2/9)&(.866, .074)&(1, .028)&&&324\\
7&4&0&(.309, .265)&(.809, .175)&(1, .06)&&&400\\
&5&2/9&(.5, .194)&(.866, .175)&(1, .06)&&&2304\\
&6&0&(.309, 1/5)&(.809, 1/5)&(1, 1/10)&&&256\\
&7&1/6&(.5, 1/6)&(.866, 1/6)&(1, 1/12)&&&1024\\
\hline
&1&1&&&&&&1\\
&2&0&(.223, .412)&(.623, .052)&(.901, .025)&(1, .010)&&49\\
&3&.352&(.5, 2/9)&(.866, .074)&(1, .028)&&&324\\
&4&0&(.223, .248)&(.623, .147)&(.901, .075)&(1, .031)&&3136\\
8&5&2/9&(.5, .194)&(.866, .139)&(1, 1/18)&&&2304\\
&6&0&(.223, .180)&(.623, .152)&(.901, .117)&(1, .051)&&12544\\
&7&1/6&(.5, 1/6)&(.866, 1/6)&(1, 1/12)&&&1024\\
&8&0&(.223, .143)&(.623, .143)&(.901, .143)&(1, .071)&&4096\\
\hline
&1&1&&&&&&1\\
&2&0&(.223, .412)&(.623, .052)&(.901, .025)&(1, .010)&&49\\
&3&.344&(.383, .213)&(.707, .063)&(.924, .037)&(1, .016)&&1024\\
&4&0&(.223, .248)&(.623, .147)&(.901, .075)&(1, .031)&&3136\\
9&5&.213&(.383, .178)&(.707, .113)&(.924, .072)&(1, .031)&&25600\\
&6&0&(.223, .180)&(.623, .152)&(.901, .117)&(1, .051)&&12544\\
&7&.156&(.383, .147)&(.707, 1/8)&(.924, .103)&(1, .047)&&65536\\
&8&0&(.223, .143)&(.623, .143)&(.901, .143)&(1, .071)&&4096\\
&9&1/8&(.383, 1/8)&(.707, 1/8)&(.924, 1/8)&(1, 1/16)&&16384\\
\hline
&1&1&&&&&&1\\
&2&0&(.174, .409)&(.5, .049)&(.766, .021)&(.940, .014)&(1, .006)&81\\
&3&.344&(.383, .213)&(.707, .063)&(.924, .037)&(1, .016)&&1024\\
&4&0&(.174, .241)&(.5, .137)&(.766, .062)&(.940, .042)&(1, .019)&14400\\
10&5&.213&(.383, .178)&(.707, .113)&(.924, .072)&(1, .031)&&25600\\
&6&0&(.174, .172)&(.5, .136)&(.766, .093)&(.940, .068)&(1, .031)&186624\\
&7&.156&(.383, .147)&(.707, 1/8)&(.924, .103)&(1, .047)&&65536\\
&8&0&(.174, .134)&(.5, .123)&(.766, .107)&(.940, .092)&(1, .043)&331776\\
&9&1/8&(.383, 1/8)&(.707, 1/8)&(.924, 1/8)&(1, 1/16)&&16384\\
&10&0&(.174, 1/9)&(.5, 1/9)&(.766, 1/9)&(.940, 1/9)&(1, 1/18)&65536\\ 
\hline\hline
\end{tabular}
\label{table:poly}
\end{table}

\end{document}